\documentclass[12pt]{article}
\usepackage{latexsym}
\usepackage{fullpage}
\usepackage{amssymb}
\begin{document}

\def\a{\alpha}
\def\b{\beta}
\def\ba{\backslash}
\def\bb{D_+}
\def\be{\begin{equation}}
\def\beq{\begin{eqnarray*}}
\def\cc{D_-}
\def\cs{\cdots}
\def\d{\tilde D}
\def\dr{{\partial_r}}
\def\dt{{\partial\over{\partial t}}}
\def\e{\varepsilon}
\def\ee{\end{equation}}
\def\eeq{\end{eqnarray*}}
\def\ep{\e^{123}_{\a\b\g}}
\def\f{\varphi}
\def\g{\gamma}
\def\gg{\bar g}
\def\im{\Longrightarrow}
\def\j{{\it{\bf j}}}
\def\ji{\,{_{^{\_\!\_}}}\!{^{_{_|}}}\,}
\def\l{\lambda}
\def\la{\langle}
\def\ld{\ldots}
\def\ms{\medskip}
\def\n{\nabla}
\def\ni{\noindent}
\def\nn{\bar{\nabla}}
\def\nt{\tilde{\nabla}}
\def\o{\omega}
\def\p{\psi}
\def\pa{\partial}
\def\pcl{\Cl (M)}
\def\SO{{\rm SO}}
\def\Hol{{\rm Hol}}
\def\hol{{\mathfrak{hol}}}
\def\SU{{\rm SU}}
\def\Spin{{\rm Spin}}
\def\Sp{{\rm Sp}}
\def\Id{{\rm Id}}
\def\U{{\rm U}}
\def\pso{P_{\SO_n}}
\def\psp{P_{\Spin_n}}
\def\psu{P_{\U(1)}M}
\def\psou{P_{\SO(n+1)}M}
\def\psoub{P_{\SO(n+1)}\bar M}
\def\psub{P_{\U(1)}\bar M}
\def\r{\rightline{$\Box$}\ss}
\def\ra{\rangle}
\def\so{\hbox{SO}}
\def\sp{\hbox{Spin}}
\def\spin{\hbox{Spin}}
\def\ss{\smallskip}
\def\t{\tilde}
\def\wt{\widetilde}
\def\te{\theta}
\def\tp{\tilde\Psi}
\def\vs{\vskip .5cm}
\def\vv{$V_+\ $}
\def\x{\times}

\def\G{\Gamma}
\def\Cl{{\bf C}l}
\def\D{\tilde D}
\def\L{\Lambda}
\def\O{\Omega}
\def\P{\Psi}
\def\Ric{{\rm Ric}}
\def\X{X^{\ast}}
\def\Y{Y^{\ast}}

\def\CM{\ensuremath{\mathbb C}}
\def\RM{\ensuremath{\mathbb R}}
\def\HM{\ensuremath{\mathbb H}}
\def\ZM{\ensuremath{\mathbb Z}}

\def\RP#1{\ensuremath{\mathbb RP^{#1}}}
\def\CP#1{\ensuremath{\mathbb CP^{#1}}}
\def\HP#1{\ensuremath{\mathbb HP^{#1}}}

\def\s{\sigma}
\def\S{\Sigma}
\def\.{\cdot}

 
\parindent0em 
\parskip1ex 
 


\newtheorem{ede}{Definition}[section]

\newtheorem{epr}[ede]{Proposition}

\newtheorem{ath}[ede]{Theorem}

\newtheorem{elem}[ede]{Lemma}

\newtheorem{ere}[ede]{Remark}

\newtheorem{ecor}[ede]{Corollary}
 

\title{Parallel spinors and holonomy groups}

\author{Andrei Moroianu and Uwe Semmelmann\footnote{Supported by the
    Deutsche Forschungsgemeinschaft}}

\maketitle

\begin{abstract} In this paper we complete the classification of spin
  manifolds admitting parallel spinors, in terms of the Riemannian
  holonomy groups. More precisely, we show that on a given
$n$--dimensional Riemannian
  manifold, spin structures with parallel spinors are in one to
  one correspondence with lifts to $\spin_n$ 
of the Riemannian holonomy group, with fixed points on the spin representation
  space. In particular, we obtain the first examples of compact
  manifolds with two different spin structures carrying parallel spinors.

\end{abstract}

\section{Introduction}

The present study is motivated by two articles (\cite{wa2}, \cite{McI2}) 
which deal with the classification of non--simply connected
manifolds admitting parallel spinors. In
\cite{wa2}, Wang uses 
representation--theoretic techniques as well as  
some nice ideas due to McInnes (\cite{McI1}) in order to obtain the
complete list of the possible holonomy groups of manifolds admitting parallel
spinors (see Theorem \ref{wang1}). We shall here be concerned with the
converse question, namely: 

(Q) Does a spin manifold whose holonomy group appears in the list above 
admit a parallel spinor ?

The first natural idea that one might have is the following
(cf. \cite{McI2}): let $M$ be a  
spin manifold and let $\t M$ its universal  cover (which is
automatically spin); let $\Gamma$ be the fundamental group of $M$ and
let  $\psp {\t M}{\rightarrow}\pso {\t M}$ be the unique spin
structure of $\t M$;  then there is a natural $\Gamma$--action on the
principal bundle $\pso {\t M}$  and the lifts of this action to 
$\psp {\t M}$ are in one--to--one correspondence with the spin 
structures on $M$.
This approach seems to us quite
unappropriated in the given context since it is very difficult to have a good
control on these lifts. 
Our main idea was to remark that the question (Q) above is not well--posed. 
Let us, indeed, consider the following slight modification of it:

(Q') If $M$ is a Riemannian manifold whose holonomy group belongs to the
list above, does $M$ admit a spin structure with parallel spinors?

It turns out that the answer to this question is simply "yes",
(see Theorem \ref{main} below). The related
question of how many such spin structures may exist on a given
Riemannian manifold is also completely solved by our Theorem \ref{n}
below. In particular, we obtain the interesting result that every
Riemannian manifold with holonomy group $\SU_m\rtimes\ZM_2$
($m\equiv0(4)$), (see explicit compact examples of such manifolds in
Section 5), has exactly two different spin structures with parallel
spinors. The only question which remains open is the existence of compact 
non--simply connected manifolds 
with holonomy $\Sp_m\x\ZM_d$ ($d$ odd and dividing $m+1$).
We remark that our results correct statements of McInnes given in
\cite{McI2}  (see Sections 4 and 5 below). 

The topic of this paper was suggested to us by K. Galicki. We
acknowledge useful discussions with D. Kotschick, A. Dessai and
D. Huybrechts. The second named author would like to thank the
IHES for hospitality and support.

\section{Preliminaries}

A spin structure on an oriented Riemannian manifold $(M^n,g)$ is a
$\Spin_n$ principal bundle over $M$, together with an equivariant
2--fold covering $\pi:\psp M\to\pso M$ over the oriented orthonormal frame bundle of $M$. Spin structures exist if and only if the second Stiefel--Whitney class $w_2(M)$ vanishes. In that case, they are in one--to--one correspondence with elements of $H^1(M,\ZM_2)$. Spinors are sections of the complex vector bundle $\S M:=\psp M\x_{\rho}\S_n$ associated to the spin structure via the usual spin representation $\rho$ on $\S_n$. The Levi--Civita connection on $\pso M$ induces canonically a covariant derivative $\nabla$ acting on spinors. 

Parallel spinors are sections  $\phi$ of $\S M$ satisfying the
differential equation $\nabla \phi \equiv 0$. They obviously correspond to fixed points (in $\S_n$) of the restriction of $\rho$ to the spin holonomy group $\wt{\Hol}(M)\subset \Spin_n$. The importance of manifolds with parallel spinors comes from the fact that they are Ricci--flat:

\begin{elem} \label{pa}{\em (\cite{hi})}
The Ricci tensor of a Riemannian spin manifold admitting a parallel spinor vanishes.
\end{elem}
{\sc Proof.}
Applying twice the covariant derivative to the parallel spinor  $\phi$
gives that the curvature of the spin--connection $\nabla$ vanishes in the
direction of  $\phi$. A Clifford contraction together with the first Bianchi identity then show that $\Ric(X)\.\phi\equiv 0$ for every vector $X$, which proves the claim.

\r

We will be concerned in this paper with irreducible Riemannian
manifolds, {\it i.e.} manifolds whose holonomy representation is
irreducible. By the de Rham decomposition theorem, a manifold is
irreducible if and only if its universal cover is not a Riemannian product.
Simply connected irreducible spin manifolds carrying parallel spinors are classified by their (Riemannian) holonomy group in the following way:

\begin{ath} \label{wang} {\em (\cite{hi}, \cite{wa1})} 
Let $(M^n, \, g)$ be a simply connected irreducible spin manifold
$(n\ge2)$. Then $M$ carries a parallel spinor if and only if the
Riemannian holonomy  
group $Hol(M, g)$ is one of the following : $G_2$ $(n=7)$; $\Spin_7$ $(n=8)$; 
$\SU_m$ $(n=2m)$; $\Sp_k$ $(n=4k)$.
\end{ath}

{\sc Proof.} If $M$ carries a parallel
spinor, it cannot be locally symmetric. Indeed, $M$ is Ricci--flat by
the Lemma above, and Ricci--flat locally symmetric manifolds are
flat. This would contradict the irreducibility hypothesis. One may
thus use the Berger--Simons theorem which states that the holonomy group of
$M$ belongs to the following list: $G_2$ $(n=7)$; $\Spin_7$ $(n=8)$; 
$\SU_m$ $(n=2m)$; $\Sp_k$ $(n=4k)$; $\U_m$ $(n=2m)$; $\Sp_1\.\Sp_k$ 
$(n=4k)$; $\SO_n$.
On the other hand, 
if $M$ carries a parallel spinor then there exists a fixed point in $\S_n$
of $\widetilde{\Hol}(M)$ and hence a vector $\xi \in \S_n$ on which the Lie
algebra $\widetilde{\hol}(M)=\hol(M)$ of $\widetilde{\Hol}(M)$ acts
trivially. It is easy to see that the spin representation of the Lie algebras
of the last three groups from the Berger--Simons list has no fixed points,
thus proving the first part of the theorem (cf. \cite{wa1}). 

Conversely, suppose that $\Hol(M)$ is one of $G_2$, $\Spin_7$,
$\SU_m$ or $\Sp_k$. In particular,
it is simply connected. Let $\pi$ denote the universal covering
$\Spin_n\to \SO_n$. Since $\Hol(M)$ is simply connected, $\pi^{-1}\Hol(M)$
has two connected components, $H_0$ (containing the unit element) and
$H_1$, each of them being mapped bijectively onto $\Hol(M)$ by
$\pi$. Now, it is known that $\pi:\wt{\Hol}(M)\to \Hol(M)$ is onto
(\cite{kn}, Ch. 2, Prop. 6.1). Moreover, $\wt{\Hol}(M)$ is connected
(\cite{kn}, Ch. 2, Thm. 4.2) and contains the unit in $\Spin_n$, so
finally $\wt{\Hol}(M)=H_0$.  
The spin representation of the Lie algebra of $H_0$ acts trivially on some
vector $\xi \in \S_n$, which implies that $h(\xi)$ is constant for $h\in H_0$. 
In particular $h(\xi)=1(\xi)=\xi$ for all $h\in H_0$, and one deduces
that $\xi$ is a fixed point of the spin representation of $H_0$. 

\r

{\bf Remark.} In the first part of the proof one has to use some
representation theory in order to show that the last three groups in the
Berger--Simons list do not occur as holonomy groups of manifolds with
parallel spinors. The non--trivial part concerns only $\U_m$ and
$\Sp_k$, since the spin representation of
$\mathfrak{so}_n=\mathfrak{spin}_n$ has of course no fixed point. 
An easier argument which excludes these two groups is the remark that they do
not occur as holonomy groups of Ricci--flat manifolds (see \cite{be}).
It is natural
to ask in this context whether there exist any simply 
connected Ricci--flat manifolds with holonomy $\SO_n$. 
Our feeling is that it should be possible to construct local examples 
but it seems to be much more difficult to construct compact examples.
Related to this, it was remarked by
A.~Dessai that a compact irreducible Ricci--flat manifold with vanishing first
Pontrjagin class must have holonomy $\SO_n$ (c.f. \cite{anand}).

\section{Wang's holonomy criterion}

In this section we recall the results of Wang (cf.~\cite{wa2})
concerning the possible holonomy groups of non--simply connected,
irreducible spin manifolds with parallel spinors. By Lemma
\ref{pa}, every  
such manifold $M$ is Ricci--flat. The restricted holonomy group
$\Hol_0(M)$ is isomorphic the full holonomy group of the universal
cover $\t M$, so it belongs to the list given by Theorem \ref{wang}. 
Using the fact that $\Hol_0(M)$ is normal
in $\Hol(M)$, one can obtain the list of all possible holonomy groups of 
irreducible Ricci--flat manifolds (see \cite{wa2}). 
If $M$ is compact this list can be considerably reduced (see \cite{McI1}).

The next point is the following simple observation of Wang 
(which we state from a slightly different point of view, more 
convenient for our purposes). It gives a criterion for a subgroup of
$\SO_n$ to be the holonomy group of a $n$--dimensional manifold with
parallel spinors:  

\begin{elem}\label{crit}
Let $(M^n, \, g)$ be a spin manifold admitting a parallel
spinor. Then there exists an embedding $\phi : \Hol(M) \rightarrow
Spin_n$ such that $\pi \circ \phi = \Id_{\Hol(M)}$. Moreover, the
restriction of the spin representation to $\phi(\Hol(M))$ has a fixed
point on $\S_n$.   
\end{elem}

Finally, a case by case analysis using this criterion yields

\begin{ath} \label{wang1} {\em (\cite{wa2})} 
Let $(M^n, \, g)$ be a irreducible Riemannian
spin manifold which is not simply connected.
If  $M$  admits a non--trivial parallel spinor, then 
the full holonomy group $\Hol(M)$ belongs to the following table:\vs

\hskip 1cm\begin{tabular}{|l|l|l|l|l|}\hline
$\Hol_\circ(M)$    &  $\dim(M)$   &  $\Hol(M)$  &  $N$   &  conditions \\ 
 \hline\hline
$\SU_m$            &    $2m$    &    $\SU_m$    &  2     &  \\  \cline{3-5}
                  &  & $\SU_m \rtimes \ZM_2$ & 1    & $m \equiv
 0(4)$ \\ \hline
           &      &     $\Sp_m$    &  $m+1$ & \\  \cline{3-5}
$\Sp_m$  & $4m$  &   $\Sp_m \times \ZM_d$   & $(m+1)/d$ & $d > 1, \; d \;
 \mbox{odd}, \; d \; \mbox{divides} \; m+1$ \\ \cline{3-5}
&  & $\Sp_m \cdot \Gamma$ & see \cite{wa2} & $m \equiv 0(2)$  \\ \cline{3-5}

\hline
$\Spin_7$          &   8  &    $\Spin_7$        &  1     &   \\
\hline
$G_2$             &    7    &    $G_2$        &  1     &             \\
\hline
\end{tabular}

\begin{center} Table 1.
\end{center}

where $\Gamma $ is either $\ZM_{2d}\ (d>1)$, or an infinite subgroup of 
      $\U(1)\rtimes\ZM_2$, or a binary dihedral, tetrahedral,
      octahedral or icosahedral group. Here $N$ denotes the
dimension of the space of parallel spinors. If, moreover, $M$ is compact, then only the following possibilities may
occur:\vs

\hskip 1cm\begin{tabular}{|l|l|l|l|l|}\hline
$\Hol_\circ(M)$    &  $\dim(M)$   &  $\Hol(M)$  &  $N$   &  conditions \\ 
 \hline\hline
$\SU_m$            &    $2m$    &    $\SU_m$    &  2     & $m$ odd \\  \cline{3-5}
                  &  & $\SU_m \rtimes \ZM_2$ & 1    & $m \equiv
 0(4)$ \\ \hline
$\Sp_m$  & $4m$  &   $\Sp_m \times \ZM_d$   & $(m+1)/d$ & $d > 1, \; d \;
 \mbox{odd}, \; d \; \mbox{divides} \; m+1$ \\ 
\hline
$G_2$             &    7    &    $G_2$        &  1     &             \\
\hline
\end{tabular}

\begin{center} Table 2.
\end{center}
\end{ath}

\section{Spin structures induced by holonomy bundles}

We will now show that the algebraic restrictions on the holonomy group 
given by Wang's theorem are actually 
sufficient for the existence of a spin structure carrying parallel spinors. 
The main tool is the following converse to Lemma \ref{crit}:

\begin{elem} \label{lift} Let $M$ be a Riemannian manifold and suppose
  that there exists an embedding $\phi:\Hol(M)\to \Spin_n$ which makes 
  the diagram
$$\begin{array}{ccccc}
& &  & &\Spin_n\cr 
&&\phi&&\cr
& &\nearrow& &\downarrow\cr
&&&&\cr
 &\Hol(M)& \longrightarrow& &\SO_n\cr
\end{array}$$
commutative. Then $M$ carries a spin structure whose holonomy group is 
exactly  $\phi(\Hol(M))$, hence isomorphic to $\Hol(M)$.
\end{elem}

{\sc Proof.} Let $i$ be the inclusion of $\Hol(M)$ into $\SO_n$ and 
$\phi:\Hol(M)\to \Spin_n$ be such that $\pi\circ\phi=i$. We fix a frame 
$u\in \pso M$ and let $P\subset \pso M$ denote the holonomy bundle of 
$M$ through $u$, which is a $\Hol(M)$ principal bundle (see \cite{kn}, 
Ch.2). There is then a canonical bundle isomorphism $P\x_i \SO_n\simeq 
\pso M$ and it is clear that $P\x_\phi \Spin_n$ together with the canonical 
projection onto $P\x_i \SO_n$ defines a spin structure on $M$. The 
spin connection comes of course from the restriction to $P$ of the 
Levi--Civita connection of $M$ and hence the spin holonomy group is 
just $\phi(\Hol(M))$, as claimed.
\r

Now, recall that Table 1  was obtained in the following way: among 
all possible holonomy groups of non--simply connected irreducible 
Ricci--flat Riemannian manifolds, one selects those whose holonomy group 
lifts isomorphically to $\Spin_n$ and such that the spin representation 
has fixed points when restricted to this lift. Using the above Lemma 
we then deduce at once the following classification result, which contains
the converse of Theorem \ref{wang1}.

\begin{ath} \label{main}
An oriented non--simply connected irreducible 
Riemannian manifold has a spin structure carrying parallel spinors if 
and only if its Riemannian holonomy group appears in Table 1 (or, 
equivalently, if it satisfies the conditions in Lemma \ref{crit}).
\end{ath}

There is still an important point to be clarified here. Let $G=\Hol(M)$ be 
the holonomy group of a manifold such that $G$ belongs to Table 1 and 
suppose that there are several lifts $\phi_i:G\to \Spin_n$ of the inclusion 
$G\to \SO_n$. By Lemma \ref{lift} each of these lifts gives rise to a spin 
structure on $M$ carrying parallel spinors, and one may legitimately ask 
whether these spin structures are equivalent or not. The answer to this 
question is given by the following (more general) result.

\begin{ath} \label{n} Let $G\subset \SO_n$ and let $P$ be a $G$--structure 
on $M$ which is connected as topological space. Then the enlargements 
to $\Spin_n$ of $P$ using two different lifts of $G$ to $\Spin_n$ are not 
equivalent as spin structures.
\end{ath}

{\sc Proof.} Recall that two spin structures $Q$ and $Q'$ are said to
be {it equivalent} if there exists a bundle isomorphism $F:Q\to Q'$
such that the diagram  

$$\begin{array}{ccccc}
Q & & \displaystyle{\mathop{\longrightarrow}^F}  & &Q'\cr 
&&&&\cr
& \searrow &&\swarrow &\cr
&&&&\cr
 & &\pso M& &\cr
\end{array}$$
commutes. 
Let $\phi_i:G\to \Spin_n$ $(i=1,2)$ be two different lifts of $G$ and
suppose that $P\x_{\phi_i}\Spin_n$ are equivalent spin structures on
$M$. Assume that there exists a bundle map $F$ which makes the diagram
$$\begin{array}{ccccc}
P\x_{\phi_1}\Spin_n & & \displaystyle{\mathop{\longrightarrow}^F}  & &P\x_{\phi_2}\Spin_n\cr 
&&&&\cr
& \searrow &&\swarrow &\cr
&&&&\cr
 & &P\x_i \SO_n& &\cr
\end{array}$$
commutative. This easily implies the existence of a smooth mapping $f:P\x \Spin_n\to \ZM_2$ such that 
\be\label{f}F(u\x_{\phi_1}a)=u\x_{\phi_2}f(u,a)a,\  \forall u\in P, \ a\in \Spin_n.\ee
 
As $P$ and $\Spin_n$ are connected we deduce that $f$ is constant, say $f\equiv \e$. Then (\ref{f}) immediately implies $\phi_1=\e \phi_2$, hence $\e=1$ since $\phi_i$ are group homomorphisms (and both map the identity in $G$ to the identity in $\Spin_n$), so $\phi_1=\phi_2$, which contradicts the hypothesis.

\r

Using the above results, we will construct in the next section the first examples of (compact) Riemannian manifolds with several spin structures carrying parallel spinors.

{\bf Remark.} Let us also note that a simple check through the list
obtained by McInnes in \cite{McI1} shows that the holonomy group of a
compact, orientable, irreducible, Ricci--flat manifold of
non--generic holonomy and real dimension not a multiple of four is
either $G_2$ or $\SU_m$ ($m$ odd) (there are two other possibilities in the
non--orientable case). Theorem 2 of \cite{McI2} (which states that
the above manifolds have a unique spin structure with parallel
spinors) follows thus immediately from our results above.

\section{Examples and further remarks}

Theorem \ref{main} is not completely satisfactory as long as 
we do not know whether for each group in Tables 1 or 2, Riemannian 
manifolds having this group as holonomy group really exist. This is 
why we will show in this section that most of the concerned groups 
have a realization as holonomy groups. We will leave as an open problem 
whether there exist compact non--simply connected manifolds 
with holonomy $\Sp_m\x\ZM_d$ ($d$ odd and $m+1$ divisible by $d$). 
We also remark that the problem which we consider here is purely 
Riemannian, {\it i.e.} does not make reference to spinors anymore.

{\bf 1. $M$ compact.} Besides the above case which we do not treat
here, it remains to construct examples of compact non--simply connected
manifolds with holonomy $G_2$, $ \SU_m$ and $\SU_m\rtimes \ZM_2$ (as these are
the only cases occurring in Wang's list in the compact case). The
first one is obtained directly using the work of Joyce (\cite{jo1}),
who has constructed several families of compact non--simply connected
manifolds with holonomy $G_2$ for which he computes explicitly the
fundamental group. 

For the second we have to find irreducible, non--simply
connected Calabi--Yau manifolds of odd complex dimension. Such examples can be 
constructed in arbitrary high dimensions. For instance, one can take the 
quotient of a hypersurface of degree $p$ in 
$\CM P^{p-1}$ by a free $\ZM_p$ action, where $p\ge5$ is any prime number 
(see \cite{beau} for details).


Finally, we use an idea of Atiyah, Hitchin and McInnes
to construct manifolds with holonomy group  $\SU_m\rtimes \ZM_2$. 
Let $a_{ij}$, ($i=1,\ld,m+1,\ j=0,\ld,2m+1$) be (strictly) positive real 
numbers and $M_i$ be the quadric in $\CP{2m+1}$ given by 
$$M_i=\{[z_0,\ld,z_{2m+1}]\ |\ \sum_{j}a_{ij}z_j^2=0\}.$$
We define $M$ to be the intersection of the $M_i$'s, and remark that
if the $a_{ij}$'s are chosen generically (i.e. such that the quadrics
are mutually transversal), then $M$ is a smooth complex
$m$--dimensional manifold realized as a complete intersection. By
Lefschetz' hyperplane Theorem (\cite{gh}) $M$ is connected and simply connected
(for $m>1$). 

Moreover, $M$ endowed with the metric inherited from $\CP{2m+1}$ becomes a 
K{\"a}hler manifold. The adjunction formula (see \cite{gh}) shows that 
$c_1(M)=0$.
Consequently $M$ is a Calabi--Yau manifold, and there exists a Ricci--flat 
K{\"a}hler metric $h$ on $M$ whose K{\"a}hler form $\O_h$ lies in the same 
cohomology class as the K{\"a}hler form $\O_g$ of $g$. We now consider the 
involution $\s$ of $M$ given by $\s([z_i])=[\bar z_i]$, which has no fixed 
points on $M$ because of the hypothesis $a_{ij}>0$. 

\begin{elem} The involution $\s$ is an anti--holomorphic isometry of $(M,h,J)$.
\end{elem}

{\sc Proof.} It is easy to see that $\s$ is actually an isometry of the Fubini--Study metric on $\CP{2m+1}$, hence $\s^*\O_g=-\O_g$. On the other hand, $\s^*h$ is a Ricci--flat K{\"a}hler metric, too, whose K{\"a}hler form is $\O_{\s^*h}=-\s^*\O_h$. At the level of cohomology classes we have thus $[\O_{\s^*h}-\O_g]=\s^*[\O_g-\O_h]=0$ and by the uniqueness of the solution to the Calabi--Yau problem, we deduce that $\s^*h=h$, as claimed.

\r
We now remark that the manifold $M$ is irreducible.  Indeed, from the
Lefschetz hyperplane theorem also follows  $b_2(M) =
1$. On the other hand, if $M$ would be reducible, the de 
Rham decomposition theorem would imply that $M=M_1\x M_2$ where $M_i$
are simply connected compact K{\"a}hler manifolds, hence $b_2(M)=b_2(M_1)+b_2(M_2)\ge2$, a
contradiction.

\begin{ecor} \label{c} The quotient $M/\s$ is a $2m$--dimensional
  Riemannian manifold  with holonomy $\SU_m\rtimes\ZM_2$.
\end{ecor}

Note that this manifold is oriented if and only if $m$ is even. For $m\equiv0(4)$, $\SU_m\rtimes\ZM_2$ has two different lifts to 
$\Spin_{2m}$, each of them satisfying the conditions of Lemma \ref{crit}. 
We thus deduce (by Theorem \ref{n}) that (for $m\equiv0(4)$) the above 
constructed manifold $M$ is a compact Riemannian manifold with exactly 
two different spin structures carrying parallel spinors.

{\bf Remark.}  This result is a counterexample to McInnes' Theorem 1 in
  \cite{McI2}, which asserts that a compact, irreducible
  Ricci--flat manifold of non--generic holonomy and real dimension
  $4m$ admits a parallel spinor if and only if it is simply
  connected. The error in McInnes' proof comes from the fact that starting
  from a parallel spinor on a manifold with local holonomy
  $\SU_{2m}$, the 'squaring' construction does not always furnish the whole
  complex volume form. In some cases one only obtain its real or complex part, 
  which is of course not sufficient to conclude that the whole holonomy 
  group is $\SU_{2m}$.

{\bf 2. $M$ non--compact.} We now give, for each group in Table 1, examples 
of (non--compact, non--simply connected) oriented Riemannian manifolds having 
this group as holonomy group. Of course, we will not consider here the 
holonomy groups of the compact manifolds constructed above, since it suffices 
to remove a point from such a manifold to obtain a non compact example. All 
our examples for the remaining groups in Table 1 will be obtained as cones 
over manifolds with special geometric structures. Recall that if $(M, \, g)$ 
is a Riemannian manifold
the cone  ${\bar M}$ is the product manifold $M \times \RM_+$ equipped with the warped product metric ${\bar g} := r^2 g \oplus d r^2$. Note that
${\bar M}$ is always a non--complete manifold, and by a result of
Gallot (\cite{gal}), the cone over a complete manifold is always 
irreducible or flat as Riemannian manifold. Using the O'Neill
formulas for warped products, it is easy to relate the different
geometries of a manifold and of its cone in the following way~(see for
example \cite{baer}, or \cite{BGM} for the definitions).

\begin{ath} {\em (\cite{baer})} Let $M$ be a Riemannian manifold and
  $\bar M$ the cone over it. Then $\bar M$ is hyperk{\"a}hler or
  has holonomy $\Spin_7$ if and only if $M$ is a 3--Sasakian manifold,
  or a weak 
  $G_2$--manifold respectively. There is an explicit natural correspondence
  between the above structures on $M$ and $\bar M$.

\end{ath}

This directly yields examples of oriented, non--simply connected
Riemannian manifolds with holonomy $\Spin_7$ and $\Sp_m$, as cones over
non--simply connected weak $G_2$--manifolds (cf.~\cite{FKMS} or
\cite{GS} for examples), and non--simply connected 3--Sasakian
manifolds respectively (cf.~\cite{BGM} for examples).

Let now $M$ be a regular simply connected 3--Sasakian manifold other
than the round sphere (all known examples of such manifolds are
homogeneous). It is a classical fact that $M$ is the total space a $\SO_3$ 
principal bundle over a quaternionic K{\"a}hler manifold, such
that the three Killing vector fields defining the 3--Sasakian structure
define a basis of the vertical fundamental vector fields of this
fibration.

For $d>1$ odd, let $\G$ be the image of $\ZM_d\subset \U(1)\subset
\SU_2$ through the natural homomorphism $\SU_2\rightarrow \SU_2/\ZM_2\simeq
\SO_3$. It is clear that $\G\simeq \ZM_d$ and by the above $\G$ acts
freely on $M$. On the other hand, for every $x\ne\pm1$ in $\SU_2$, the right
action of $x$ on $\SO_3$ preserves a one--dimensional space of left
invariant vector fields and defines a non--trivial rotation on the
remaining 2--dimensional space of left invariant vector
fields on $\SO_3$. This means that if $\g\ne1$ is an arbitrary element of $\G$,
its action on $M$ preserves exactly one Sasakian structure and defines
a rotation on the circle of Sasakian structures orthogonal to the
first one. The following classical result then shows that the holonomy group of
the cone over $M/\Gamma$ has to be $\Sp_m\x\ZM_d$.

\begin{epr}
Let $(M, \, g)$  be a Riemannian manifold with universal cover
${\t M}$. If the natural surjective homomorphism
$ \pi_1(M) \rightarrow \Hol/\Hol_0$ is not bijective, then there exists
a subgroup $K \subset \pi_1(M)$ such that ${\t M}/K$
is a manifold with  $\Hol = \Hol_0$. The group $K$ is actually the
kernel of the homomorphism above.
\end{epr}
 
Similarly one may construct examples of manifolds with holonomy
$\Sp_m\.\G$ for every group $\G$ listed in Theorem \ref{wang1}.

 \labelsep .5cm

\ni(A. Moroianu) Centre de Math{\'e}matiques, Ecole Polytechnique, UMR 7640 du CNRS, Palaiseau, France\\
\ni(U. Semmelmann) Mathematisches Institut, Universit{\"a}t M{\"u}nchen, Theresienstr. 39, 80333, M{\"u}nchen, Germany.

\end{document}